# $Z_8$ IS NOT DUALIZABLE

## CS. SZABÓ

## 1. Introduction

In [4] natural (strong) duality is proved for the ring $\mathbf{Z}_4$. They also show it for $\mathbf{Z}_{p^2}$, where $p$ is a prime. The next question is, whether $\mathbf{Z_8}$ is dualizable or not. There were several attempts to approach this problem. The closest shot was made by Lousindi Sabourin, who interpreted the problem into a question of quadratic equations over vector spaces. Let $V$ be a vector space over $F_2$, the two element field. A subset $S \subseteq V$ has the property

($\mathcal{Q}$) if $S$ is the set of solutions of a quadratic equation

($\mathcal{Q}_3$) if $S \cap W$ is the set of solution of a quadratic equation for every 3-dimensional affine subspace $W$ of $V$.

Of course, $\mathcal{Q}$ implies $\mathcal{Q}_3$. Sabourin showed ([6]) that if $\mathcal{Q}_3$ implies $\mathcal{Q}$ (Sindi's conjecture), then $\mathbf{Z_8}$ admits a natural duality.

In this paper we show that $\mathbf{Z_8}$ does not admit a natural duality. In fact, we show that $2\mathbf{Z_8} = \{2,4,6,8 \,|\, +,\cdot\}$ is not dualizable, and this will imply that the original ring is not dualizable, either. As a corollary we show that Sindi's conjecture does not hold. Our technique will be similar to the one in [5], where non-dualizability is proved for the quaternion group.

This work was delivered when the author spent a year at The Fields Institute for Research in Mathematics, in Toronto, Ontario, as a post. doc. of Matt Valeriote and Bradd Hart during the Algebraic Model Theory year. The topic was taken up by Ross Willard while making efforts to get people work on duality. The final inspiration for this work came at a conversation with Steve Seif ([7]). The ideal circumstances for work was provided by the staff of The Fields Institute and we would like to say special thanks to Pauline Grant and Becky Sappong.

## 2. Remarks on duality

This chapter is supposed to be skipped by the ones who have already experienced some duality before.

*Date*: July 24, 1997.
The research of the author was supported by a grant from the NSERC of Canada.

For the benefit of readers not familiar with the theory of natural dualities, we begin with a brief review of what is meant by *'admitting a (natural) duality'* and refer to Davey [2] or the forthcoming text Clark and Davey [1] for a detailed account.

Let $\mathbf{A}$ be a finite algebra and let $\tilde{A} = \langle A; F, H, R, \tau \rangle$ be a topological structure on the same underlying set $A$, where:

(a) each $f \in F$ is a homomorphism $f : \mathbf{A}^n \to \mathbf{A}$ for some $n \in \mathbb{N} \cup \{0\}$,
(b) each $h \in H$ is a homomorphism $h : dom(h) \to \mathbf{A}$ where $dom(h)$ is a subalgebra of $\mathbf{A}^n$ for some $n \in \mathbb{N}$,
(c) each $r \in R$ is (the universe of) a subalgebra of $\mathbf{A}^n$ for some $n \in \mathbb{N}$,
(d) $\tau$ is the discrete topology.

Whenever (a), (b) and (c) hold, we say that the operations in $F$, the partial operations in $H$ and the relations in $R$ are *algebraic over* $\mathbf{A}$. These compatibility conditions between the structure on $\mathbf{A}$ and the structure on $\tilde{A}$ guarantee that there is a naturally defined dual adjunction between the quasivariety $\mathcal{A} := \mathbb{ISP}\mathbf{A}$ generated by $\mathbf{A}$ and the topological quasivariety $\mathcal{X}_{\tilde{A}} := \mathbb{IS}_c\mathbb{P}\tilde{A}$ generated by $\tilde{A}$; if there is no chance of confusion, we will write $\mathcal{X}$ for $\mathcal{X}_{\tilde{A}}$. For all $\mathbf{B} \in \mathcal{A}$ the homset $D(\mathbf{A}) := \mathcal{A}(\mathbf{B}, \mathbf{A})$ of all homomorphisms from $\mathbf{B}$ to $\mathbf{A}$ is a closed substructure of the direct power $\tilde{A}^B$ and for all $\tilde{X} \in \mathcal{X}$ the homset $E(\tilde{X}) := \mathcal{X}(\tilde{X}, \tilde{A})$ is a subalgebra of the direct power $\mathbf{A}^X$. It follows easily that the contravariant hom-functors $\mathcal{A}(-, \mathbf{A}) : \mathcal{A} \to \mathcal{S}$ and $\mathcal{X}(-, \mathfrak{A}) : \mathcal{X} \to \mathcal{S}$, where $\mathcal{S}$ is the category of sets, lift to contravariant functors $D : \mathcal{A} \to \mathcal{X}$ and $E : \mathcal{X} \to \mathcal{A}$.

For each $\mathbf{B} \in \mathcal{A}$ there is a natural embedding $e_{\mathbf{B}}$ of $\mathbf{B}$ into $ED(\mathbf{B})$ given by evaluation: for each $b \in B$ and each $x \in D(\mathbf{B}) = \mathcal{A}(\mathbf{B}, \mathbf{A})$ define $e_{\mathbf{B}}(b)(x) := x(b)$. Similarly, for each $\tilde{X} \in \mathcal{X}$ there is an embedding $\varepsilon_{\tilde{X}}$ of $\tilde{X}$ into $DE(\tilde{X})$. A simple calculation shows that $e : \text{id}_{\mathcal{A}} \to DE$ and $\epsilon : \text{id}_{\mathcal{X}} \to DE$ are natural transformations. If $e_{\mathbf{B}}$ is an isomorphism for all $\mathbf{B} \in \mathcal{A}$ we say that $\tilde{A}$ yields a *(natural) duality* on $\mathcal{A}$. If there is some choice of $F$, $H$ and $R$ such that $\tilde{A}$ yields a duality on $\mathcal{A}$ then we say that $\tilde{A}$ (or $\mathcal{A}$) *admits a natural duality* or, briefly, is *dualizable*.

## 3. Acceccories

We wish to prove that for no choice of $F$, $H$ and $R$ does $\tilde{Z}_8$ yield a duality on $\mathcal{G}$, the quasivariety generated by $\mathbf{Z_8}$. For this, it is enough to show that there is no duality when $F = H = \emptyset$ and $R$ consists of all subgroups of all finite powers of $\mathbf{Z_8}$, the so-called *brute force duality*, see [1] or [2]. In order to prove that there is no brute force duality, we need to find a (necessarily infinite) group $\mathbf{D} \in \mathcal{G}$ such that $e_{\mathbf{D}}$ is not onto $ED(\mathbf{D})$. We will use the *ghost element* method. We will choose $\mathbf{D}$ to be a subring of $\mathbf{Z_8}^{\mathbb{Z}}$ and will define a continuous structure preserving map $\Phi : D(\mathbf{D}) \to \tilde{Z}_8$ such that $\{\Phi(\pi_i) \mid i \in \mathbb{Z}\}$ is not an element of $\mathbf{D}$, implying that $\Phi$ is not the evaluation map for any element of $\mathbf{D}$ and therefore that $e_{\mathbf{D}}$ is not onto $ED(\mathbf{D})$. Here $\pi_i$ is the (restriction to $\mathbf{D}$ of the) i-th projection of $\mathbf{Z_8}^Z$ onto $\mathbf{Z_8}$.

For a subring $\mathbf{D}$ of $R^Z$ let $\mathbf{D}_{fin}$ denote the elements of $\mathbf{D}$ with finitely many nonzero coordinates. We say that an element $\bar{\mathbf{v}}$ of $\mathbf{D}$ has *eventual value* $v$ (is eventually $v$) if all its coordinates but finitely many ones are equal to $v$. So $\mathbf{D}_{fin}$ is the set (subring) of eventually 0 elements of $\mathbf{D}$.

3.1. **The ring.** The ring $\mathbf{Z_8}$ is an algebra with two binary operations and a constant.

$$\mathbf{Z_8} = \{1, 2, 3, 4, 5, 6, 7, 0 \mid +, \cdot, 0\}.$$

The Jacobson(nil)-radical of $\mathbf{Z_8}$,

$$J(\mathbf{Z_8}) = 2\mathbf{Z_8} = \{0, 2, 4, 6 \mid +, \cdot, 0\},$$

is a two-class nilpotent ring.

First we construct a subring $\mathbf{D}$ of $2\mathbf{Z_8}^Z$ with the ghost-vector missing, and after that we examine the possible homomorphisms from $\mathbf{D}$ to $\mathbf{Z_8}$.

The vectors $\bar{b} = (2, 2, 0)$, $\bar{a} = (0, 2, 2)$ generate a subring $R$ of $\mathbf{Z_8}^3$, isomorphic to the 2 generated free ring in the variety generated by $J(\mathbf{Z_8})$. Define $\bar{\mathbf{b}}$ and $\bar{\mathbf{a}}_i$ for $i \in Z$, elements of $R^Z$ as follows:

$$\bar{\mathbf{b}}_i = \bar{b} \quad \text{and} \quad \bar{\mathbf{a}}_{ij} = \begin{cases} \bar{b} & \text{if } i = j \\ \bar{0} & \text{if } |i - j| = 1 \\ \bar{a} & \text{if } |i - j| > 1 \end{cases}$$

Let $\mathbf{D}_1 = \langle \bar{\mathbf{b}}, \bar{\mathbf{a}}_i \mid i \in Z \rangle$. Moreover, let $\bar{\mathbf{e}}_i = \bar{\mathbf{a}}_i - \bar{\mathbf{a}}_{i-1}$, then

$$\bar{\mathbf{e}}_{ij} = \begin{cases} \bar{a} & \text{if } j = i - 2 \\ -\bar{b} & \text{if } j = i - 1 \\ \bar{b} & \text{if } i = j \\ -\bar{a} & \text{if } j = i + 1 \\ \bar{0} & \text{otherwise,} \end{cases}$$

and so, $\mathbf{D}_1 = \langle \bar{\mathbf{b}}, \bar{\mathbf{a}}_0, \bar{\mathbf{e}}_i \mid i \in \mathbf{Z} \rangle$, i.e., $\mathbf{D}_1$ is generated by the vectors of the form:

$\bar{\mathbf{b}} = (\ldots, \bar{b}, \bar{b}, \bar{b}, \ldots),$

$\bar{\mathbf{a}}_0 = (\ldots, \bar{a}, \bar{a}, \bar{a}, 0, \overset{0}{\bar{b}}, 0, \bar{a}, \bar{a}, \ldots),$

$$\bar{\mathbf{e}}_i = (\ldots, 0, 0, \bar{a}, -\bar{b}, \overset{i}{\bar{b}}, , -\bar{a}, 0, 0, \ldots).$$

Let $\mathbf{D}_2 = \langle \bar{\mathbf{x}}^2 \,|\, \bar{\mathbf{x}} \in R^Z \rangle_{fin}$, the ring generated by the squares of the elements in $R^Z$ containing finitely many nonzero coordinates. Observe that $\mathbf{D}_2 = \langle \bar{a}^2, \bar{b}^2 \rangle^{\mathbf{Z}}_{fin}$, and $\mathbf{D}_2 \cdot R^Z = 0$ holds. Finally, let $\mathbf{D} = \langle \mathbf{D}_1, \mathbf{D}_2 \rangle = \mathbf{D}_1 + \mathbf{D}_2$.

### 3.2. The ghost-vector.

Our ghost vector will be $\overline{\mathbf{ab}}$, where $\overline{\mathbf{ab}}_i = \bar{a} \cdot \bar{b} = \bar{a}\bar{b}$. First we have to show that $\overline{\mathbf{ab}}$ is not in $\mathbf{D}$, but it can be 'approximated' by vectors from $\mathbf{D}$.

**Claim** The vector $\overline{\mathbf{ab}} = (\ldots, \bar{a}\bar{b}, \bar{a}\bar{b}, \bar{a}\bar{b}, \ldots)$ is not contained in $\mathbf{D}$, but $_i\overline{\mathbf{ab}} = (\ldots, \bar{a}\bar{b}, \overset{i}{\bar{0}}, , \bar{a}\bar{b}, \ldots)$, where $\bar{0}$ is at the $i$-th place is in $\mathbf{D}$.

*Proof of the claim.* The variety generated by $2\mathbf{Z_8}$ satisfies the identities $2x = x^2$ and $4x = (= 2x^2) = 0$. Thus a generator set of $\langle \bar{a}\bar{b} \rangle^Z \cap \mathbf{D}$ can be obtained on the following way: Take the product of all pairs of generators and substitute each coordinate equal to $\bar{a}^2$ or $\bar{b}^2$ with 0. All these elements satisfy the following 'parity check condition': Every vector is eventually 0 or eventually $\bar{a}\bar{b}$. In the first case the sum of the coordinates is 0 (there are even many $\bar{a}\bar{b}$-s), in the second case there are odd many 0-s. Obviously, this property is preserved by addition of these elements and so, the property holds for $\langle \bar{a}\bar{b} \rangle^Z \cap \mathbf{D}$, proving the claim. □

As $R$ has a natural embedding into $\mathbf{Z_8}$, $R^Z$ has a natural embedding to $(\mathbf{Z_8}^3)^Z$, that gives a natural embedding of $\mathbf{D}$ into $R^Z$. For this subring and for its elements we shall use the notations above, we shall denote the elements of $\mathbf{D}$ and their images at the embedding on the same way. (e.g. $\bar{\mathbf{b}} = (\ldots, 0, 2, 2, 0, 2, 2, \ldots)$ and $\overline{\mathbf{ab}} = (\ldots, 0, 4, 0, 0, 4, 0, \ldots)$. Thus $\mathbf{D} \leq \mathbf{Z_8}^Z$.

### 3.3. The maps.

First we examine the possible maps from $\mathbf{D}$ to $\mathbf{Z_8}$. Let $f \in \hom(\mathbf{D}, \mathbf{Z_8})$. Since $\mathbf{D}$ satisfies the identity $4x = 0$, the same holds for its image, hence $\mathbf{D}$ is mapped to $2\mathbf{Z_8}$. We are interested in the action of $f$ on the set $\{\,_i\overline{\mathbf{ab}}\,|\,i \in Z\,\}$, in fact, we will show that $f(_i\overline{\mathbf{ab}}) = f(_j\overline{\mathbf{ab}})$ for almost all $i, j \in Z$. Note that $_i\overline{\mathbf{ab}} - _j\overline{\mathbf{ab}} \in \mathbf{D}_{fin}$. Moreover, $\bar{\mathbf{e}}_i \cdot \bar{\mathbf{e}}_{i+2} = (\ldots, 0, 0, 0, \overset{i}{\bar{a}\bar{b}}, \overset{i+1}{\bar{a}\bar{b}}, 0, 0, \ldots)$ and so, $\langle \bar{a}\bar{b} \rangle^Z \cap \mathbf{D}_{fin} = \langle \bar{\mathbf{e}}_i \cdot \bar{\mathbf{e}}_{i+2} \,|\, i \in Z \rangle$.

**First case:** the image of $\mathbf{D}_{fin}$ at $f$ is a zeroring. Then $f(x)f(y) = f(xy) = 0$ for any $x, y \in \mathbf{D}_{fin}$, so $f(\bar{\mathbf{e}}_i \cdot \bar{\mathbf{e}}_{i+2}) = 0$. Thus $f(\langle \bar{a}\bar{b} \rangle^Z \cap \mathbf{D}_{fin}) = \{0\}$ holds, hence $f(_i\overline{\mathbf{ab}}) = f(_j\overline{\mathbf{ab}})$ for every $i, j \in Z$.

**Last case:** the image of $\mathbf{D}_{fin}$ is not a zeroring. Then there is an $i \in Z$ such that $f(\bar{\mathbf{e}}_i) = 2$ or $6$. As $\bar{\mathbf{e}}_i \cdot \bar{\mathbf{e}}_j = 0$ if $|i-j| > 3$, $f(\bar{\mathbf{e}}_j)$ is contained in $\{0,4\}$, the annihilator of $2\mathbf{Z_8}$ for $|i-j| > 3$. Thus there is a smallest $i$ such that $f(\bar{\mathbf{e}}_i) = 2$ or $6$, and for $i < j$, $f(\bar{\mathbf{e}}_j)\mathbf{Z_8} = 0$. Without loss of generality we may assume that $f(\bar{\mathbf{e}}_1) = 2$ but $f(\bar{\mathbf{e}}_i) \in \{0,4\}$ for $i < 1$. As $\bar{\mathbf{e}}_1\bar{\mathbf{e}}_j = 0$ for $j \geq 3$, $f(\bar{\mathbf{e}}_j)$ annihilates $2\mathbf{Z_8}$ in this case, too. From this easily can be derived that for $\bar{\mathbf{h}} \in \langle \bar{a}\bar{b}\rangle^Z \cap \mathbf{D}_{fin}$, if $\bar{\mathbf{h}}_i = 0$ for $i < 2$ or $\bar{\mathbf{h}}_i = 0$ for $i > 2$, then $f(\bar{\mathbf{h}}) = 0$. Thus $f(_i\overline{\mathbf{ab}}) = f(_j\overline{\mathbf{ab}})$ if $2$ is not between $i$ and $j$.

Notice that $f(\bar{\mathbf{e}}_0)f(\bar{\mathbf{e}}_1) + f(\bar{\mathbf{e}}_0)f(\bar{\mathbf{e}}_2) = 0$ as each summand contains the image of a vector with index smaller than $1$. On the other hand $f(\bar{\mathbf{e}}_0)f(\bar{\mathbf{e}}_1) + f(\bar{\mathbf{e}}_0)f(\bar{\mathbf{e}}_2) = f(\bar{\mathbf{e}}_1)f(\bar{\mathbf{e}}_2 + \bar{\mathbf{b}})$. Since $f(\bar{\mathbf{e}}_1)$ generates $2\mathbf{Z_8}$, the element $f(\bar{\mathbf{e}}_2) + f(\bar{\mathbf{b}}) \in Ann(\mathbf{Z_8})$, and so $f(\bar{\mathbf{e}}_2)(f(\bar{\mathbf{e}}_2) + f(\bar{\mathbf{b}})) = 0$. Also, by examining the indices, $f(\bar{\mathbf{e}}_2)f(\bar{\mathbf{e}}_{-1}) + f(\bar{\mathbf{e}}_2)f(\bar{\mathbf{e}}_5) = 0$. Adding the last two equalities we get $f(\bar{\mathbf{e}}_2)(f(\bar{\mathbf{e}}_2) + f(\bar{\mathbf{e}}_{-1}) + f(\bar{\mathbf{e}}_5) + f(\bar{\mathbf{b}})) = f(\ldots, 0, 0, 0, \overset{0}{\bar{a}\bar{b}}, 0, 0\overset{3}{\bar{a}\bar{b}}, 0, 0, 0, \ldots) = 0$. This is a vector with two $\bar{a}\bar{b}$ entries on different sides of the critical 2nd coordinate, and so if non of $i$ and $j$ equals $2$, $f(_i\overline{\mathbf{ab}}) = f(_j\overline{\mathbf{ab}})$. We showed that for each $f$, which $f(\mathbf{D}_{fin}) = 2\mathbf{Z_8}$ holds for, there is a coordinate $i(f)$, such that $f(_l\overline{\mathbf{ab}}) = f(_m\overline{\mathbf{ab}})$ holds if $m$ and $l$ are different from $i(f)$. We call $i(f)$ the *critical coordinate* of the map $f$.

## 4. The results

Now, we know everything to prove our main theorem:

**Theorem 4.1.** *The ring $\mathbf{Z_8}$ does not admit a natural duality*

*Proof.* Let $\phi$ be the following map form $\hom(\mathbf{D},\mathbf{Z_8})$:

$$\phi(f) = \begin{cases} f(_0\overline{\mathbf{ab}}) & \text{if } f \text{ belongs to the first case,} \\ f(_{i+1}\overline{\mathbf{ab}}) & \text{if } i \text{ is the critical coordinate of } f. \end{cases}$$

In order to show that $\phi$ is structure preserving, for any finite set of homomorphisms $f_1, \ldots, f_n$ from $\mathbf{D}$ to $\mathbf{Z_8}$ we have to find an element $\bar{\mathbf{v}}$ of $\mathbf{D}$ such that $\phi(f_i) = f_i(\bar{\mathbf{v}})$ for $i = 1, \ldots, n$. Let $m$ be different (if any exists) from the critical coordinates of the $f_i$-s.

Then by Section 3.3 $\phi(f_i) = f(_m\overline{\mathbf{ab}})$ for all these maps.

For showing that $\phi$ is continuous, it is enough (and easy) to see that if two map agrees on $_0\overline{\mathbf{ab}}$, $_1\overline{\mathbf{ab}}$ and $_2\overline{\mathbf{ab}}$, then they have the same value at $\phi$.

On the other hand, $\phi(\pi_i) = \bar{a}\bar{b}$ and $\overline{\mathbf{ab}}$ is not in $\mathbf{D}$, hence $\phi$ is not an evaluation map. □

**Corollary 4.2.** Sindi's conjecture fails. There is a vector space $V$ over $F_2$ and a subset $S$ of $V$, such that for every 3-dimensional subspace $W$ of $V$, the intersection $W \cap S$ is the solution set of some quadratic equation, but there is no quadratic equation with $S$ as its the solution set.

Several times rings are considered as algebras with a unit element 1. In that case
$$\mathbf{Z_8} = \{1, 2, 3, 4, 5, 6, 7, 0 \,|\, +, \cdot, 0, 1\,\}$$
an algebra with two constants, 0 and 1. For us, the main difference is that in this case $\mathbf{D}$ has to contain the $\bar{\mathbf{1}} = (\ldots, 1, 1, 1, \ldots)$ constant 1 vector. If we try to add it to our construction, the ring turns out to contain $\overline{\mathbf{ab}}$, the ghost vector. But, if at the beginning we take $\bar{a} = (2, 2, 0, 0,)$ and $\bar{b} = (0, 2, 2, 0)$, the whole idea of the construction can be saved, and the proof goes exactly the same way.

**Corollary 4.3.** $\mathbf{Z_8}$ with 1 is not dualizable.

FIELDS INSTITUTE/ MCMASTER UNIVERSITY
*E-mail address*: `csaba@@cs.elte.hu`